\documentclass[11pt,a4paper,leqno]{amsart}

\usepackage{fullpage,amsmath,amsthm,amsfonts,amscd,
latexsym,amssymb,eucal,curves}
  \usepackage[all]{xy}
\usepackage[dvips]{graphicx}
\usepackage[dvips]{color}
\newcommand{\ZZ}{\mathbb{Z}}

\newcommand{\CC}{\mathbb{C}}

\newcommand{\PP}{\mathbb{P}}
\newcommand{\NN}{\mathbb{N}}

\newcommand{\HH}{\mathbb{H}}

\newcommand{\QQ}{\mathbb{Q}}
\newcommand{\RR}{\mathbb{R}}

\newcommand{\PicS}{{\rm Pic}}

\newcommand{\Jac}{{\rm Jac} \,}
\newcommand{\cLL}{{\mathcal L}}

\newcommand{\SL}{{\rm SL}}

\newcommand{\diag}{{\rm diag}}

\newcommand{\norm}{{\rm norm}}
\newcommand{\sms}{\smallsetminus}

\newcommand{\ol}{\overline}

\newcommand{\GL}{{\rm GL}}

\newtheorem{Defi}{Definition}[section]

\newtheorem{Rem}[Defi]{Remark}

\newtheorem{Prop}[Defi]{Proposition}
\newtheorem{Lemma}[Defi]{Lemma}
\newtheorem{Cor}[Defi]{Corollary}

\newtheorem{Example}[Defi]{Example}
\newtheorem{Thm}[Defi]{Theorem}

\begin{document}{\large}
\title{Finiteness results for Teichm\"uller curves}
\date{\today}
\author{Martin M\"oller}

\begin{abstract} 
We show that for each genus there are only finitely
many algebraically primitive Teichm\"uller curves $C$,
such that i) $C$ lies in the hyperelliptic locus
and ii) $C$ is generated by an abelian differential with 
two zeros of order $g-1$.
\newline
We prove moreover that for these Teichm\"uller curves 
the trace field of the affine group is not only totally
real but cyclotomic.
\end{abstract}
\thanks{Supported by the 
DFG-Schwerpunkt ``Komplexe Mannigfaltigkeiten''}
\maketitle

\section*{Introduction}

A Teichm\"uller curve is an algebraic curve in the moduli
space of curves of genus $g$, denoted by $M_g$, whose
preimage in Teichm\"uller space is
a complex geodesic for the Teichm\"uller metric. Teichm\"uller
geodesics are obtained as the orbit of a pair $(X^0,q^0)$ of 
a Riemann surface $X^0$ plus a quadratic differential $q^0$ 
on $X^0$ under the
action of $\SL_2(\RR)$. Those (few) pairs $(X^0,q^0)$ that give
Teichm\"uller curves are called Veech surfaces. We restrict 
ourselves to the case when $q^0 = (\omega^0)^2$ is a square 
of a holomorphic differential. The case
of proper squares might be analysed using the canonical double
covering of $X^0$, that makes the pullback into a square.
We remark that Teichm\"uller curves naturally lift to
the bundle $\Omega M_g$ over $M_g$ of holomorphic one-forms. This bundle
is stratified according to the multiplicity
of the zeros of the one-form. See Section \ref{Not} for more details.
\par
The first examples of Teichm\"uller curves were obtained as
coverings of the torus ramified over one point. There are 
infinitely many of them in each stratum. First examples 
not of this type were discovered by Veech (\cite{Ve89}).
In particular the trace field of the affine group (see Section
\ref{Not}) is not $\QQ$ in Veech's examples.
\newline
In genus $2$ there are infinitely many non-torus coverings
in the stratum with one double zero but only a single one in the stratum
with a two zeros. If one fixes one additional discrete parameter
(the discriminant of the order all curves parametrized by 
such a Teichm\"uller curve have real multiplication with),
the number becomes finite also for Teichm\"uller curves in 
the stratum with one double zero. In fact there are one or two of 
them according to the congruence class of the discriminant mod $8$.
This classification is contained in \cite{McM04a}, 
\cite{McM04b} and \cite{McM04c}.
\par
To go beyond genus $2$ we recall from \cite{Mo04a} that the
family of Jacobians over a Teich\-m\"uller curve splits 
into an $r$-dimensional part with real multiplication and
some rest, where $r$ is the field extension degree of the
trace field over $\QQ$. 
\newline
A Teichm\"uller curve in $M_g$ is called {\em algebraically primitive}
if the trace field has degree $g$ over $\QQ$. This implies
that the curve is {\em geometrically primitive}, i.e.\ that
the pair $(X^0,\omega^0)$ does not arise from a surface
of lower genus plus a differential via a covering construction.
Both notions coincide in genus two, but in general the converse 
implication is not true.
\par
At the time of writing the following is known about
primitive Teichm\"uller curves:
\newline
Among Veech's examples there are infinitely many algebraically
primitive ones, but at most one for each genus. Besides this
there are series of examples and sporadic ones in 
\cite{Vo96}, \cite{Wa98}, \cite{KeSm00}. Only finitely many of them
are algebraically primitive and for each genus there are
only finitely many examples. The recent work of McMullen
(\cite{McM05}) contains infinitely many geometrically primitve (though
not algebraically primitive examples) for the genera $3$,
$4$ and $5$.
\par
The purpose of the present work is to obtain some finiteness
results valid in all genera. We cannot hope for such results for
imprimitive curves without fixing additional discrete parameters.
For geometrically primitive but algebraically imprimitive Teichm\"uller
curves it seems unclear what to expect. If we restrict to algebraically
primitive Teichm\"uller curves we show, generalizing \cite{McM04c}:
\par
{\bf Theorem \ref{finiteTHM}}: For fixed genus $g$ there are only 
finitely many algebraically primitive
Teichm\"uller curves in the connected component of the
stratum $\Omega M_g(g-1,g-1)$, that parametrizes hyperelliptic curves.
\par 
We will consider the family of curves $f: X \to C$ over a Teichm\"uller
curve $C$ or over a suitable cover of $C$. Recall from \cite{Mo04b} that
the zeros of the generating differential $\omega^0$ determine sections
of $f$. In the algebraically primitive case the difference of any two
of those sections is a torsion element of the relative Jacobian.
The theorem is an instance of the philosophy that torsion points
on families of curves are rare. It might be possible to show the
same type of result for differentials with more zeros instead of 
hyperelliptic ones. But using the same methods the combinatorics
become quite complicated then.
\par
We briefly outline the strategy of our proof: \par
\begin{itemize}
\item[i)] From an argument in \cite{Mo04b} we deduce that an 
algebraically primitive
Teich\-m\"uller curve in $\Omega M_g(g-1,g-1)^{\rm hyp}$
has a reducible and an irreducible degeneration (Thm.\ \ref{sepperpts})
in say the vertical and horizontal direction.
\item[ii)] The irreducible degeneration is used to bound the torsion
order (Prop.\ \ref{Nbounded} and Section \ref{Nboundedproof}). This limits the
suitably normalized {\em widths} of the cylinders in the horizontal 
direction to a finite set. It generalizes
the discussion of sine ratios in Section $2$ of \cite{McM04c}. 
Prop.\ \ref{Nbounded} has the flavour of the toric case
of the Mordell-Lang conjecture. Yet none of the versions in
the literature seems strong enough to cover what we need.
\item[iii)] The reducible degneration is used to relate the
torsion order and the {\em  moduli} of the cylinders 
in the vertical direction (Thm.\ \ref{divN}).
\item[iv)] The combination of these informations limits the possibilties
for the flat geometry of a Veech surface to a finite number
(see the prototype in Figure $2$ and the end of Section \ref{finsec}).
\end{itemize}
\par
As a byproduct of the proof we obtain:
\par
{\bf Corollary \ref{tracecor}:} The trace field of an algebraically primitive
Teichm\"uller curve in the stratum $\Omega M_g(g-1,g-1)^{\rm hyp}$ 
is cyclotomic.
\par
The cyclotomic fields appear roughly as follows: The normalisations
of some degenerate fibres in the family over the Teichm\"uller curve
are isomorphic to $\PP^1$. Arranging the position of the zeros of
the generating differential suitably, the preimages of the nodes
are forced by the torsion condition to lie at roots
of unities in $\PP^1$. We deduce that enough periods of $\omega^0$
lie in this cyclotomic field to conclude that the trace field is
cyclotomic. 
\par
We remark that the trace fields of all presently known
Teichm\"uller curves are cyclotomic. Based on the above Corollary
one might conjecture that this holds in general, at least for
Teichm\"uller curves with more than one zero.
\par
The author thanks F.~Herrlich for valuable
comments and C.~McMullen for pointing to the work of H.~Mann
on relations of roots of unities and for the useful
Lemma \ref{resLemma}.
\par

\section{Notation} \label{Not}

{\bf Strata of $\Omega M_g$}
\newline
We denote the tautological bundle over $M_g$ by $\Omega M_g$.
Its points are pairs $(X^0, \omega^0)$ of a Riemann surface
$X^0$ of genus $g$ and a holomorphic differential (or
equivalently: a one-form) $\omega^0 \in \Gamma(X^0,\Omega^1_{X^0})$.  
This space is naturally
stratified by the type of multiplicities of the zeros of
$\omega$. Kontsevich and Zorich have determined the connected
components of the strata (\cite{KoZo03}). 
\newline
A pair $(X,\omega)$ belongs to a {\em hyperelliptic stratum } 
if $X$ is a hyperelliptic curve with involution $\sigma$ and quotient map 
$\pi: X \to X/\langle \sigma \rangle \cong \PP^1$, such that
$\omega^2 = \pi^* q$ for a quadratic differential $q$ on $\PP^1$
with (a) $2g+1$ simple poles and a zero of order $2g-3$
or (b) $2g+2$ simple poles and a zero of order $2g-2$.
In case (a) the pair belongs to $\Omega M_g(2g-2)$ while
in case (b) the pair belongs to $\Omega M_g(g-1,g-1)$.
The hyperellipic strata form connected components (\cite{KoZo03} Thm.~1).
They will be denoted by a superscript $hyp$. For
$(X,\omega) \in \Omega M_g(g-1,g-1)^{{\rm hyp}}$ the involution
$\sigma$ interchanges the two zeroes of $\omega$. 
\par
Note that there are other types of zeros of a pair $(X,\omega)$
such that $X$ is hyperelliptic and such that $\omega^2$ is the
pullback of a quadratic differential on the quotient. But the
above two cases are the only ones, where a connected component consists
entirely of such hyperelliptic pairs.
\newline
The two hyperelliptic strata are the natural generalisation of
the only two strata that exist for $g=2$.
\par
{\bf $\SL_2(\RR)$-action}
\newline
There is a natural action of $\SL_2(\RR)$ on $\Omega M_g$ minus
the zero section: Apply the $\RR$-linear transformation to the local 
complex charts of $X$ given by integrating $\omega$ to obtain
a new complex structure and apply the $\RR$-linear transformation
to the real and imaginary parts of $\omega$ to obtain a new one-form,
which is holomorphic for the new complex structure. For more
details see e.g.\ \cite{MaTa02} or 
\cite{McM03}. This action obviously preserves the
stratification of $\Omega M_g$.
\par
{\bf Teichm\"uller curves} 
\newline 
Teichm\"uller curves are algebraic curves $C \to M_g$ in the moduli
space of curves of that
are geodesic for the Teichm\"uller metric. We deal here exclusively
with Teichm\"uller curves generated by a pair $(X^0,\omega^0)$
i.e.\ whose natural lift to the bundle of quadratic differentials
over $M_g$ lies in the image of $\Omega M_g$. Here $C = \HH/\Gamma$,
where $\Gamma \subset {\rm PSL}_2(\RR)$ is the image of the affine group
of $(X^0, \omega^0)$ (see e.g.\ \cite{Mo04a}).
Let $K = \QQ({\rm tr}(\gamma), \gamma
\in \Gamma)$ be the trace field of $\Gamma$ and $r:= [K:\QQ]$.
Let $f: X \to C$
denote the universal family over some finite unramified
cover of $C$, abusively denoted by the same letter. Let
$\Jac(f): \Jac X/C \to C$ denote the family of Jacobians. Recall
from \cite{Mo04a} that $\Jac X/C$ splits up to isogeny into a
product of a family $g:A \to C$ of abelian varieties of dimension $r$
with real multiplication by $K$ and a family of abelian varieties
of dimension $g-r$. Since the splitting up to isogeny is not unique
we take $g: A \to C$ to be the maximal quotient in its isogeny class.
This letter should cause no confusion with the genus of $X^0$.
\par
We extend all the above families to families over $\ol{C}$,
i.e.\ let $\ol{f}: \ol{X} \to \ol{C}$ be the stable model and
$\widetilde{f}: \widetilde{X} \to \ol{C}$ the minimal semistable model
with smooth total space $\tilde{X}$. Also let 
$\ol{g}: \ol{A} \to \ol{C}$ be the corresponding family
of semiabelian varieties.
\par 
{\bf N\'eron Models for families of Jacobians}
\newline
Let $F$ denote the function field of the curve $C$. The 
N\'eron Model $\tilde{g}: Q \to \ol{C}$ of a family 
$\ol{g}: \ol{A} \to \ol{C}$ of semiabelian varieties is a group
scheme, whose fibre over the generic point of $\ol{C}$
coincides with $\ol{A}_F$ and
such that for any given smooth group scheme $\ol{Y} \to
\ol{C}$ any map $\ol{Y}_F \to Q_F$ over $F$ 
extends uniquely to a map $\ol{Y} \to Q$ over all $\ol{C}$.
In particular sections of
$g$ extend to sections of $\tilde{g}$.
In all the cases we consider N\'eron Models exist, see \cite{BLR90}.
\newline 
In case of algebraically primitive Teichm\"uller curves, i.e.\ for
$g=r$, the family $\ol{g}$ is just 
$\PicS^0(\ol{X}/\ol{C})$, i.e.\ line bundles on $\ol{X}$ that
are of degree zero on each component of each fibre. The connected
component of $1$ of $Q$, denoted by $Q^0$, coincides with
$\PicS^0(\ol{X}/\ol{C})$ in this case (\cite{BLR90} Thm.~9.5.4 b)).
\par
{\bf Torsion}
\newline
Let $f: X \to C$ be the universal family over a Teichm\"uller
curve generated by $(X^0, \omega^0)$ 
in the stratum $\Omega M_g(k_1,\ldots,k_r)$. Recall
from \cite{Mo04b} that, maybe after passing to a
finite unramified cover of $C$, the zeros of $\omega^0$ define sections
$s_1,\ldots,s_r$ of $f$. For any pair $(i,j)$ the 
difference $s_i -s_j$ is a 
torsion section of $g$. It extends to a section of $\tilde{g}$.
Since (in characteristic zero) the kernel of multiplication 
by some integer is \'etale on any group scheme, in
particular on the N\'eron Model (\cite{BLR90} Lemma 7.3.2), the order of 
$(s_i - s_j)$ restricted to any fibre of 
$\tilde{g}$ equals the same number $N(i,j)$. 
In particular this holds for the fibres over the cusps.
\par

\section{Degenerations} \label{degensection}
We study the degenerate fibres and give a relation between
the geometry of a degenerate fibre and the
torsion order of the difference of the two zeroes,
if $(X^0,\omega^0) \in \Omega M_g(g-1,g-1)$ generates
an algebraically primitive Teichm\"uller curve.
\par  
\begin{Thm} \label{deggenus}
Let $\ol{f}: \ol{X} \to \ol{C}$ be the universal
family over a Teichm\"uller curve.The sum of the genera of the components
of a singular fibre of $\ol{f}$ is at most $g-r$. In particular
the degenerate fibres of an algebraically primitive
Teichm\"uller curve have only rational
components.
\end{Thm}
\par
{\bf Proof:} A family of abelian varieties with real 
multiplication degenerates to a semi-abelian variety 
whose abelian part is trivial (see e.g. \cite{Go02}
Lemma 2.23). Hence the abelian part of the fibre of $\ol{g}$
over any cusp has dimension at most $g-r$.
\newline
Alternatively this can be deduced from the
Clemens-Schmid exact sequence for a degeneration of
Hodge structures and the explicit description of the
local system in \cite{Mo04a}.
\hfill $\Box$
\par
We recall how the degeneration of a Teichm\"uller
curve is described via the euclidian geometry 
defined by $(X^0, \omega^0)$: A geodesic on $X^0$
has a well-defined slope and all geodesics with this slope
form a {\em  direction}.
Veech dichotomy (\cite{Ve89}) states that each direction that contains
a geodesic joining two zeros or one zero to itself (a {\em saddle
connection}) is {\em periodic} i.e.\ each geodesic in this direction
is periodic or closed.
\newline
The closed geodesics of a periodic direction (say the
horizontal one) sweep out cylinders
$C_i$ and we denote their core curves by $\gamma_i$. Consider
the degenerate fibre obtained by applying $\diag(e^t,e^{-t})$
to $(X^0, \omega^0)$ for $t \to \infty$. Say this point corresponds
to the cusp $c \in \ol{C} \sms C$. By \cite{Ma75} the
stable model of the singular (or 'degenerate') 
fibre $X_c$ of $f$ is obtained by squeezing
the core curves of the $C_i$ to points. Topologically the 
irreducible components of $X_c$ are obtained by 
cutting along the $\gamma_i$.
\par
\begin{Cor} \label{noofcyl}
Each direction of a Teichm\"uller curve
in $\Omega M_g(k_1,\ldots,k_s)$ has at least $r$ and
at most $r+s-1$ cylinders.
\end{Cor}
\par
{\bf Proof:}
Each component of the degeneration in the
given direction contains at least one zero. \hfill $\Box$
\par
For the rest of this section we suppose that $(X^0, \omega^0)$
generates {\em an algebraically primitive Teichm\"uller curve.}
\par
\begin{Thm} \label{sepperpts}
For any two zeros $Z_1$ and $Z_2$ of $\omega^0$ with
$Z_1 \neq Z_2$ there is a direction
with the following property: 
\newline
Let $X_c$ denote the singular fibre corresponding to 
the degeneration in this direction
and $s_i$ the sections defined by the $Z_i$. Then $s_1$ and
$s_2$ intersect $X_c$ in different irreducible components.
\end{Thm}
\par
{\bf Proof:}  
We know that $s_i$ does not intersect the degenerate fibre
in a node. Suppose the statement was wrong.
Then $s_1 - s_2$ 
defines a non-zero section of $\PicS^0(\ol{X}/\ol{C})$ over
the completed Teichm\"uller curve $\ol{C}$. This is not
possible by \cite{Mo04b} Thm.~3.1: Its proof shows not
only that there are only finitely many sections of $g$,
but also that there are none of $\ol{g}$.
\hfill $\Box$
\par
Suppose from now on that the Teichm\"uller curve $C$ is generated
by a differential with two zeros of order $g-1$.
By Thm.~\ref{sepperpts} there is a direction, say the
vertical one, such that the corresponding singular 
fibre $X_v$ has two components. The vertical direction has
hence $g+1$ cylinders. Let $\gamma_i$ denote the core curves
of the cylinders. We number them in such a way that
for $i=1,\ldots,a$ the curve $\gamma_i$ degenerates to a
node on the first component of $X_v$, while for $i=a+1,\ldots,a+b$
the curve $\gamma_i$ degenerates to a node on the
second component. We enumerate
the components of $X_v$ such that $a \leq b$. Note that
$a+b \leq g-1$ since the two components of $X_v$ intersect
in at least two points: a core curve of a cylinder is not separating.
\par
We denote by $h^v_i$ the height and by $b^v_i$ the width of
the $i$-th vertical cylinder, i.e.~the length of $\gamma_i$. 
Moreover let $m^v_i = h^v_i/b^v_i$ be the modulus of the
$i$-th vertical cylinder.
\newline
It is remarked in  Veech \cite{Ve89} that
the moduli $m^v_i$ for $i=1,\ldots,g$ are commensurable. It 
is no loss of generality
for the purposes below to rescale the generating differential
of the Teichm\"uller curve such that $m^v_i \in \NN$ and
${\rm gcd}(m^v_i, i=1,\ldots,g)=1$. 
\newline
A small simple loop in $C$ around the cusp $c$ obtained by
degenerating in the vertical direction
corresponds (compare \cite{Ve89} Prop.\ 2.4) to the product
$$ (\prod_{i=1}^{g+1} D_{\gamma_i}^{m^v_i})^k,
\quad  \text{where} \; 
D_{\gamma_i}\;\text {is a Dehn twist along}\;\gamma_i. $$
Here $k$ is some positive integer, which appears since we have
taken (with abuse of notation) coverings of the 
Teich\-m\"uller curve $C$ that may ramify
at the cusps. Hence the loop is not necessarily a generator
of the corresponding parabolic subgroup of the affine group.
This means that in the stable model the node in $X_v$ 
corresponding to $\gamma_i$ is given by 
$xy = t^{m^v_i k}$, where $t$ is a local coordinate
of $C$ at the cusp $c$ and $x,y$ are local coordinates of an
embedding of a neighborhood of the node in the stable fibre into $\CC^3$. 
In fact, the statement is local in the base $C$ and in the total
space and it reduces for $m^v_i k = 1$ to the easiest case of
the Picard-Lefschetz transformation. The general case is obtained
via base change. After resolving this singularity the fibre 
$\widetilde{X_v}$ of the semistable model with smooth
total space of $\ol{f}$ contains a chain of $m^v_i k -1 $ rational
$(-2)$-curves in the preimage of the node.
\par
\begin{Thm} \label{divN}
Let $N$ denote the order of $s_2 - s_1$. Suppose that the 
moduli $m^v_i$ are integers and
let $m^v = {\rm lcm} \{m^v_i, i= a+b+1,\ldots,g+1\}$.
Then the torsion order and the moduli of the cylinders
are related by 
$$ \sum_{i= a+b+1}^{g+1} \frac{m^v}{m^v_i} \,\; \text{\rm divides} \; N.$$
\end{Thm}
\par
{\bf Proof:} By the preceding discussion the fibre $\widetilde{X_v}$
looks as in the following figure. Lines correspond to components
of the semistable fibre, intersection points are nodes and $Z_1$, $Z_2$
are the intersection points of the section $s_i$ with the stable
fibre.
\begin{figure}[h] \label{semistab}
\centerline{\includegraphics{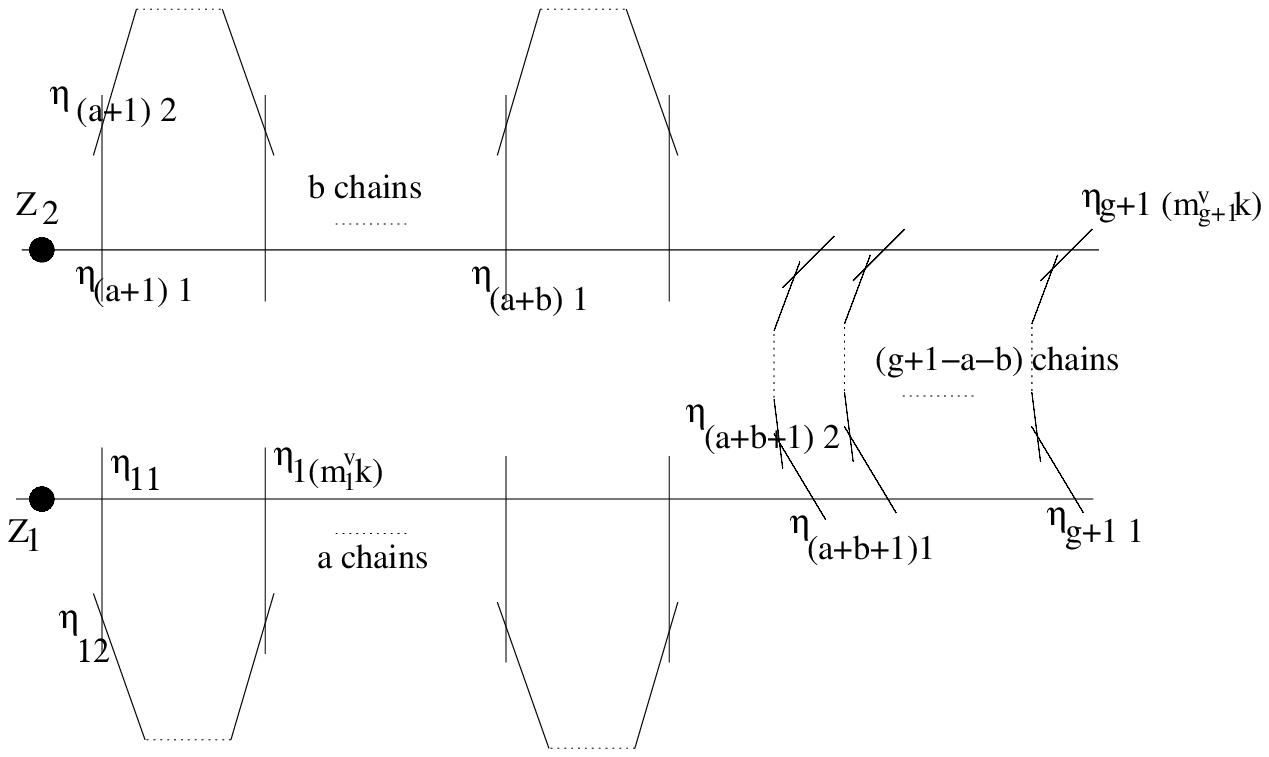}}
\caption{Semistable model of $\widetilde{X_v}$}
\end{figure}
Section 9.6 in \cite{BLR90}, in particular p.~283, gives a presentation
of the component group $Q/Q^0$ of the N\'eron model of a flat family
of curves with smooth total space over a discrete valuation ring. 
This can be applied to the localization of $\tilde{f}$ at the 'vertical' cusp: 
Let $\eta_{ij}$ denote the nodes of the singular fibre, i.e.\
the edges of the dual (intersection) graph used in loc.\ cit.
The component group is generated by the $\eta_{i1}$ and 
$\eta_{ij} - \eta_{i\; j-1}$ for $i=1,\ldots,g+1$ and
$j=2,\ldots, m^v_i k$ with the following relations:
The components of $\widetilde{X_v}$ contribute
\begin{eqnarray*}
\eta_{ij} - \eta_{i\; j-1} &=& 0 \quad\quad\quad  
i=1,\ldots,g+1,\; j=2,\ldots, m^v_i k \\
\sum_{i=a+b+1}^{g+1} \eta_{i1} + 
\sum_{i=1}^{a} (\eta_{i1} - \eta_{i\;m^v_ik}) &=& 0 \\
-\sum_{i=a+b+1}^{g+1} \eta_{i\;m^v_ik} + 
\sum_{i=a+1}^{a+b} (\eta_{i1} - \eta_{i\;m^v_ik}) &=& 0 \\
\end{eqnarray*}
and the fundamental group of the intersection graph contributes
the relations
\begin{eqnarray*}
\sum_{j=1}^{m^v_i\,k} \eta_{ij} &=& 0 \quad \quad \quad 
i=1,\ldots, a+b \\
\sum_{j=1}^{m^v_i\, k} \eta_{ij} - 
\sum_{j=1}^{m^v_{g+1} k} \eta_{g+1\,j} & = &0 
\quad \quad \quad 
 i=a+b+1,\ldots,g \\
\end{eqnarray*} 
The difference $s_2 -s_1$ defines a section of $Q$, hence of $G:=Q/Q^0$, 
which is given in this presentation e.g.\ by 
$$[s_2 -s_1 ] = \sum_{j=1}^{m^v_{g+1}\, k} \eta_{g+1\,j} $$
\par
We shall show that the order of $[s_2 -s_1]$ in $G$ equals
$\sum_{i= a+b+1}^{g+1} \frac{m^v}{m^v_i}$.
\newline
We may simplify the presentation of $G$ using only
the generators $\eta_{i,1}$ for $1,\ldots,g+1$ and relations
$$ 
\begin{array}{lcll}
(m^v_i k)\, \eta_{i,1} &=& 0 & \quad \quad \quad i=1,\ldots, a+b \\
(m^v_i k)\, \eta_{i,1} - (m^v_g k) \eta_{g+1,1} &=& 0 & 
\quad \quad \quad i=a+b+1,\ldots, g+1 \\
\sum_{i=a+b+1}^{g+1} \eta_{i,1} &=& 0 & \\ \end{array} $$
In this presentation $[s_2-s_1]= (m^v_i k)\, \eta_{i,1}$ for
any $i=a+b+1,\ldots, g+1$.
We have 
$$(\sum_{i= a+b+1}^{g+1} \frac{m^v}{m^v_i})[s_2-s_1]
=  (m^v k) \sum_{i= a+b+1}^{g+1} \eta_{i,1} =0.$$
To see that the order is not smaller, we consider $[s_2-s_1]$
in the group $H$ with the same generators and all but the
last relation. If $n\cdot[s_2-s_1] = 0$ in $G$ then there
is $n' \in \ZZ$ such that $n\cdot[s_2-s_1] = n' \sum_{i=a+b+1}^{g+1} 
\eta_{i,1}$ in $H$.
Listing the equivalence class of the right hand side in $H$ this 
means that we can write $n = \sum_{i=a+b+1}^{g+1} n_i$ such that
$$n\cdot[s_2-s_1] = \sum_{i=a+b+1}^{g+1} n_i (m^v_i k)\, \eta_{i,1} = 
\sum_{i=a+b+1}^{g+1} n' \eta_{i,1}. $$
Hence $n'$ is a common multiple of all $m_i^v k$ and hence $m^v k$, which
appeared above, is minimal. 
\par
Since the order of $s_2-s_1$ in the component group divides $N$
we are done. \hfill $\Box$
\par
\begin{Example} {\rm In case of the decagon, 
the unique primitive Teichm\"uller curve in
$\Omega M_2(1,1)$ one has $N= 5$ and the moduli are $(1,2,1)$ 
(see \cite{McM04c}). This is confirmed by $$(2/1 + 2/2 + 2/1) \;| \;5.$$}
\end{Example}
\section{Algebraically primitive
Teichm\"uller curves in $\Omega M_g(g-1,g-1)^{\rm hyp}$}
\label{finsec}
In this section we prove the following theorem:
\par
\begin{Thm} \label{finiteTHM}
There are only finitely many algebraically primitive 
Teichm\"uller curves in the
component $\Omega M_g(g-1,g-1)^{\rm hyp}$  for each $g\geq 2$. 
\end{Thm}
\par
We specialize the results of Section \ref{degensection} to the
algebraically primitive case and the hyperelliptic stratum.
We start with a direction that contains a saddle connection
joining the two zeros, say the horizontal one. 
By Thm.~\ref{deggenus} and Cor.~\ref{noofcyl}
this direction contains precisely $g$ cylinders.
Similarly as for the vertical cylinders we denote 
by $h^h_i$, $b^h_i$ and $m^h_i$ 
the (respective) heights, widths and moduli of the
horizontal cylinders.
\par
Suppose that the vertical direction is chosen as in the paragraph
preceding Thm.~\ref{divN}. 
Since the hyperelliptic involution interchanges
the zeros it also interchanges the components of such
a degenerate fibre and hence $a=b$. Moreover since
the hyperelliptic involution  on a smooth hyperelliptic curve has $2g+2$
fixed points, these fixed points have to degenerate to
$g+1$ nodes joining the two components. We have shown:
\par
\begin{Lemma} For a degenerate fibre of an
algebraically primitive Teichm\"uller curve
in the component $\Omega M_g(g-1,g-1)^{\rm hyp}$ we have
$a=b=0$.
\end{Lemma}
\par
We describe a prototype for a Veech surface in
$\Omega M_g(g-1,g-1)^{\rm hyp}$:
Suppose the hyperelliptic involution fixes $n$ of the
$g$ horizontal cylinders and interchanges the remaining
$g-n$ cylinders in pairs. This implies that we have
precisely $2n$ Weierstra\ss \ points contained in
the interiors of the horizontal cylinders. The Weierstra\ss \
points on the boundary define sections of the family $f$
(again maybe after passing to an unramified cover of $C$)
that do not pass through the nodes of the degenerate fibre.
Hence they are fixed points of the 'hyperelliptic'
involution that acts on the normalization of the degenerate fibre.
Since this normalization is isomorphic to $\PP^1$ there
are precisely two fixed points and hence precisely two
Weierstra\ss \ points on the boundary. To obtain
$2g+2$ Weierstra\ss \ points alltogether we must have $n=g$,
i.e.\ all the horizontal cylinders are fixed by the involution.
We conclude that such a Veech surface looks as in Figure $2$.
The dots correspond to the Weierstra\ss\ points, the
square and the cross denote the zeroes of $\omega^0$.
Vertical edges are glued by horizontal translations. The
horizontal edges containing the Weierstra\ss\ points are glued
on the same horizontal cylinder. In the other cases the 'free' top
horizontal saddle connection of the $i$-th cylinder is glued
to the 'free' bottom saddle connection of the $i+1$-st cylinder.
For $g$ even the square and the star have to be switched in
the lowest parallelogram.
\par
\begin{figure}[h] \label{prototype}
\centerline{\includegraphics{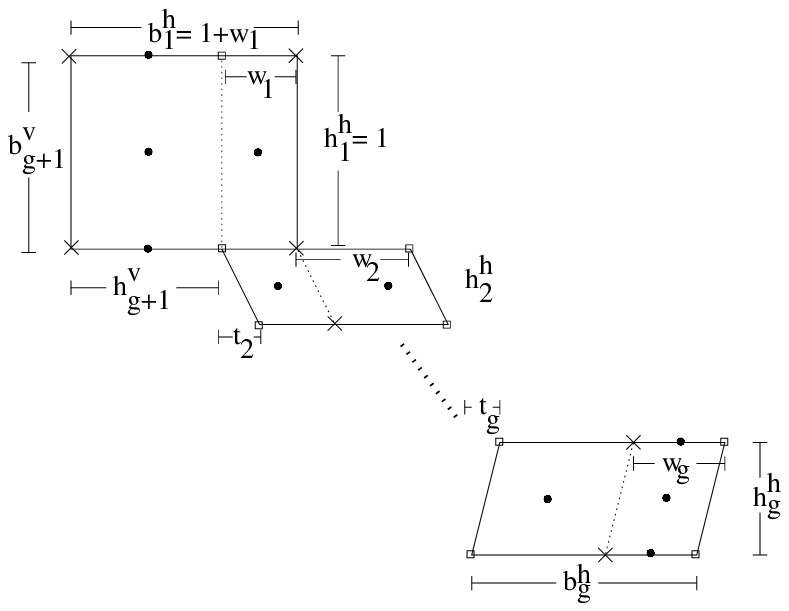}}
\caption{Prototype of a Veech surface in 
$\Omega M_g(g-1,g-1)^{\rm hyp}$}
\end{figure}
We normalize the prototype (by $\GL_2^+(\RR)$-action) by imposing that
\begin{equation} \label{Norm}
h^h_1 = 1 \quad \text{and} \quad b^h_1 = 1 + w_1, \quad \text{i.e.\ that}
\; m^v_{g+1} = 1, 
\end{equation}
where the 
the $g+1$-st vertical cylinder sits in the upper left corner of the
picture. We suppose no longer 
(as we did in Section \ref{degensection}) that the vertical
moduli are all integers  but only that they are
rational.
\par
\begin{Rem}{\rm By \cite{Mo04a} Thm.\ 2.6 the Jacobian of $X^0$
has real multiplication and the prototype above looks similar to
the ones in \cite{McM04a}, Section $3$. Nevertheless we do not claim that
the elliptic curves obtained by glueing the horizontal slits
are isogenous (what is known to be true in genus $2$).
}\end{Rem}
\par
{\bf Horizontal degeneration.}
\newline
Let $X_h$ be the singular fibre obtained by degeneration in the
horizontal direction. By Thm.~\ref{deggenus} the degenerate fibre
$X_h$ is a singular rational curve with $g+1$ nodes. We denote
by $X_h^\norm$ its normalization. We suppose
that the intersection of $X_h$ with $s_2$ and $s_1$ 
lift to the points $0$ and $\infty$ on $X_h^\norm$, respectively.
Since the hyperelliptic involution interchanges the two sections, 
we may suppose that one of its fixed points is $1$, i.e.\
that it acts by $z \mapsto 1/z$. The Weierstra\ss\ points degenerate
to $\pm 1$ and $g$ pairs $(x_i,1/x_i)$ on $X_h^\norm$ that are glued 
together on $X_h$.
\newline
Let $\cLL \subset f_* \omega_{\ol{X}/\ol{C}}$ 
be the distinguished subbundle, whose restriction to $X^0$ is
$\CC \cdot \omega^0$. We choose a generator $\omega_h$  of $\cLL|_{X_h}$
and denote
by $\omega_h^\norm$ its pullback to $X_h^\norm$. The differential
$\omega_h^\norm$  has zeros of order $g-1$ at $0$ and infinity
and simple poles at $x_i^{\pm 1}$ such that the residues differ
by a factor $-1$.
\par
Since $\Jac(f) = \ol{g}$ for an algebraically primitive curve, 
the torsion condition implies that there is a map
$t: X \to \PP^1$ whose fibre over $0$ (resp.\ $\infty$)
equals $s_1$ (resp.\ $s_2$) with multiplicity $N$. The map 
extends to $t: X_{\rm bl} \to \PP^1$ for some suitable blowup
of $\ol{X}$. Since $X_h$ is irreducible and contains
both a point that maps to $0$ and to $\infty$, the map $t$
has to be non-constant on $X_h$. Hence $t$ has the form
$z \mapsto z^N$ on $X_h^\norm$ and factors through $X_h$.
This implies that $x_i$ are $2N$-th roots of unity.
\par
Since the degenerate fibre can be obtained from $(X^0, \omega^0)$
by applying $\diag(e^t, e^{-t}) \subset \SL_2(\RR)$
the residues of $\omega^\norm_h$ around the $g$ poles coincide
up to a common scalar multiple with the integrals of $\omega^0$
along the core curves of the horizontal
cylinders, i.e.\ with the $b_i^h$.
\par
\begin{Lemma} \label{resLemma}
The residues $b^h_i$ ($i=1,\ldots,g$) of $\omega_h^\norm$, normalized
such that $b^h_1=1$, form a basis of $K/\QQ$.
\end{Lemma}
\par
{\bf Proof:} Let $\lambda$ denote a primitive element of $K/\QQ$ and
consider it (\cite{Mo04a} Thm.\ 2.6) as an endomorphism $T_\lambda$
of the family of semiabelian
varieties $\ol{g}$. In particular $T_\lambda$ acts on $X_h$. 
The differential $\omega_h$ is an eigenform for the 
action of $T_\lambda$. Hence $T_\lambda$ acts $\QQ$-linearly on 
the periods $b^h_i$ of $\omega_h$ and
$(b^h_1,\ldots,b^h_g)$ form an eigenvector for the eigenvalue 
$\lambda \in K$.
Since $b^h_1=1$ we can express all powers of $\lambda$ as $\QQ$-linear
combinations of the $b^h_i$. Since $\lambda$ is primitive and
$[K:\QQ]=g$ the $b^h_i$ form a basis. \hfill $\Box$
\par
We have shown that we can write $\omega_h^\norm$ in two ways
 $$ \omega_h^\norm =\sum_{i=1}^g \left( \frac{b^h_i}{z-x_i}
 - \frac{b^h_i}{z-x_i^{-1}} \right)dz
= C \frac{z^{g-1}}{\prod_{i=1}^{g} \left((z-x_i)(z-x_i^{-1})\right)} dz,$$
where $x_i$ are roots of unity, $b_i^h$ form a basis of a real
subfield $K \subset \QQ(x_1,\ldots,x_g)$ and $C$ is some
real number.
\par
\begin{Prop} \label{Nbounded}
For fixed $g$, there are only finitely many $g$-tupels
$(x_1,\ldots,x_g)$ of roots of unity such that there exist
real numbers $b_1^h,\ldots,b_g^h$ which form a $\QQ$-basis
of some real number field $K \subset \QQ(x_1,\ldots,x_g)$
and a real constant $C$, such that we have the following
identity of rational functions
$$ \sum_{i=1}^g \left( \frac{b^h_i}{z-x_i}  - \frac{b^h_i}{z-x_i^{-1}}
\right)
= C \frac{z^{g-1}}{\prod_{i=1}^{g} \left((z-x_i)(z-x_i^{-1})\right)}$$
In particular the least common multiple of the orders of the $x_i$
satisfying the above condition ist bounded by a function depending
only on $g$.
\end{Prop}
\par
The proof does not require any properties of Teichm\"uller curves
and will be given in the next section. 
\par
\begin{Cor} \label{widthcor}
There is only a finite number 
of period tupels $(b^h_2,\ldots,b^h_g,w_1)$
that can occur for a Veech curve normalized as in (\ref{Norm}).
\end{Cor}
\par
{\bf Proof:} The finiteness of possibilities for the $b^h_i$
is an immediate consequence.
The period $w_1$ is the integral of $\omega_h$ along a path
from $0$ to $\infty$ that crosses the unit circle once
(in a point different from $\pm 1$ and $x_i^{\pm 1}$). 
Since the $x_i$ are fixed up to a finite number of choices
this determines $w_1$ up to a finite number of choices.
\hfill $\Box$
\par
{\bf Vertical direction}
\newline
The work has been done in the previous section. We record that
Thm.~\ref{divN} implies:
\par
\begin{Lemma} \label{modulifinite}
For fixed $N$ there is only a finite number of possibilities
for the moduli $m^v_i$.
\end{Lemma}
\par
{\bf Proof of Thm.~\ref{finiteTHM}:} 
Fix one of the finitely many possibilities for
the $b^h_i$, $w_1$ and hence $w_i$ ($i=1,\ldots,g$) and
for the moduli $m^v_j$ ($j=1,\ldots,g$).
For all $j$ the heights $h^v_j$ are bounded above
by ${\rm max}\{w_i; i = 1,\ldots,g-1\}$. 
Hence all the $b^v_j$ are
bounded above.
\newline
Let $J_1 \subset \{1,\ldots,g\}$ be the indices of vertical
cylinders intersecting $w_1$. 
For $j \in J_1$ the heights $h^v_j$ are bounded away from zero
since $b^v_j$ is bounded away from zero by $h^h_1=1$ and the $m^v_j$ are
fixed. Since the $b^h_i$ are fixed and $h^v_i$ are bounded
away from zero, there is a only finite number of possibilities for the 
intersection numbers $e_{ij} := \gamma^h_i \cdot \gamma^v_j$ for 
$i = 1,\ldots,g$ and $j \in J_1$. We fix one possibility.
\par
{\em Claim:} For at least one (say the $i_0$-th) of 
the horizontal cylinders intersected by some $j \in J_1$ 
the height $h^h_{i_0}$ is bounded away from zero by a constant
depending only on $w_1$, the moduli $m^v_j$ and the
intersection numbers fixed so far.
\newline
In fact, we know that 
$$ w_1 = \sum_{j \in J_1} e_{1j} h^v_j$$ and by definition
$$ h^v_j = m^v_j \sum_{i=1}^{g} e_{ij} h^h_j.$$
Putting these equations together we obtain using $h^h_1=1$
$$ w_1 - \sum_{j \in J_1} m^v_j e^2_{1j}
= \sum_{i=2}^{g} h^h_i \cdot (\sum_{j \in J_1}  m^v_i e_{1j} e_{ij}).$$
The left hand side of this equation is non-negative and if
it were zero this would imply that the vertical cylinders crossing
$w_1$ do not intersect any other horizontal cylinder but the first.
This is absurd. Hence it is positive and depends only on quantities
fixed so far. This implies that not all the $h^h_i$ for
$i=2,\ldots,g$ can be simultaneously arbitrarily small.
\par
Using the claim we denote by $J_2$ the set of
cylinders that intersect $\gamma_1$ or $\gamma_{i_0}$. 
As above this limits the $e_{ij}$ for $j \in J_2$ to a 
finite number. We now proceed inductively analysing 
$w_j$ in place of $w_1$ to conclude that all 
intersection numbers $e_{ij}$ vary through a finite list.
\par 
Fix one of the finitely many possibilities for the
intersection numbers. We know that 
$$ b^h_i = \sum_{j=1}^{g+1} e_{ij} h^v_j =
\sum_{j=1}^{g+1} e_{ij} m^v_j b^v_j$$
for $i=1,\ldots,g$ and for $j=1,\ldots,g$ we know by definition 
$$ b^v_j = \sum_{i=1}^{g} e_{ij} h^h_i.$$
\par
Let $E$ denote the $g\times (g+1)$-matrix with entries $e_{ij}$.
From \cite{HuLa05} we deduce that $E \diag(m^v_i) E^t$ is regular.
In fact they show that the eigenvalues of $E \diag(m^v_i) E^t \diag(m^h_i)$
form a basis of $K/\QQ$. 
Hence we may plug the second equation above in the first and
solve uniquely for the $h^h_i$, since we know the $b^h_i$.
This also determines the $b^v_j$ and consequently the $h^v_j$.
\par
We know all heights and widths of the cylinders and it remains to
limit the possible twists $t_i$ for $i=2,\ldots,g$ to a finite number.
The absolute values of the twists are bounded by the intersection
numbers times $b^h_i$ and they can only vary in positive integral 
linear combinations of the $h^v_i$. Hence there is only a finite number
of possibilities for the twists.
\hfill $\Box$
\par
\begin{Cor} \label{tracecor}
The trace field $K$ of an algebraically primitive Teichm\"uller
curve in $\Omega M_g(g-1,g-1)^{\rm hyp}$ is abelian.
\end{Cor}
\par
{\bf Proof:} In the above proof of Thm.\ref{finiteTHM}
we have seen that the periods $b^h_i$ lie in the field $\QQ(x_i)$,
where $x_i$ are roots of unity. 
The field generated by the $b^h_i$ coincides with the
trace field of $\Gamma$ by Lemma \ref{resLemma}.
\hfill $\Box$
\par
Since it fits into this context, we end our finiteness
discussion by the following complement:
\par
\begin{Thm} Fix a genus $g$ and consider all Teichm\"uller
curves $C \to M_g$. If we fix moreover the Euler characteristic
$\chi(C) = 2g-2+n$ of $C$ then there is only a finite number
of possibilities for the monodromy, in particular for
the trace field of such a Teichm\"uller curve.
\end{Thm}
\par
{\bf Proof:} The Euler characteristic of the corresponding
curve in the moduli space $M_g^{[3]}$ of curves with level-$3$-structure
is also bounded. Hence we can apply Prop.\ 3.10 in \cite{De87}.
\hfill $\Box$
\par
\section{Proof of Prop.\ \ref{Nbounded}} \label{Nboundedproof}
Suppose we are given a rational function as in the statement of
Prop.\ \ref{Nbounded}. Choosing $\prod_{i=1}^g (z-x_i)(z-x_i^{-1})$
as common denominator and comparing coefficients of $z^0$ to
$z^{g-2}$ translates
into the following system of equations for $e=1,\ldots,g-1$
(the coefficients of $z^g$ to $z^{2g-2}$ provide the
same system):
\begin{equation}
\sum_{i=1}^g \left(b_i^h (x_i-x_i^{-1}) \sum_{
\stackrel{j_1<\ldots<j_e}{{\rm all}\,j_k \neq i}} \;
\prod_{k=1}^{e-1} (x_{j_k}+x_{j_k}^{-1}) \right)= 0
\end{equation}
We subtract in the first step
$\sum_{j=1}^g (x_j+x_j^{-1})$ times the equation with $e=1$ from 
$e=2$ to obtain an equation denoted by $(Eq:2')$. Then subtracting 
$\sum_{j_1<j_2} (x_j+x_j^{-1})$ times the equation with $e=1$ from 
$e=3$ and adding $\sum_{j=1}^g (x_j+x_j^{-1})$ times the equation
$(Eq:2')$ we obtain an equation denoted by $(Eq:3')$. Proceeding in this
way we obtain the simplified system
$$ \sum_{i=1}^g \left(b_i^h (x_i-x_i^{-1})(x_i+x_i^{-1})^{e-1}
\right) = 0. \eqno(Eq:e')$$
This system of equations is equivalent to the system
$$\sum_{i=1}^g \left(b_i^h (x_i^e-x_i^{-e})
\right) = 0. \eqno(Eq:e'') $$
for $e=1,\ldots,g-1$, which will be used in the sequel.
\par
We say that an equation 
\begin{equation} \label{principalrel}
\sum_{i=1}^{k} a_i \zeta_i = 0 
\end{equation}
where the $\zeta_i$ are pairwise different roots of unity and where 
the $a_i$ lie in the number field $K$ form a {\em $K$-relation
of length $k$}. The relation is called irreducible, if 
$\sum_{i=1}^k b_i \zeta_i =0$ and $b_i(a_i-b_i)=0$ for
all $i$ implies that $b_i=0$ for all $i$ or $a_i-b_i=0$ for all $i$.
Each relation is a sum of irreducible relations, but there may be
several ways of writing a relation as sum of irreducible
relations.
\par
\begin{Lemma} \label{MannThm}
Let $\sum_{i=1}^k a_i \zeta_i=0$ be an irreducible $K$-relation
with $K \subset \QQ^{\rm ab}$ and $[K:\QQ] = g$. Then multiplying the
relation by a suitable root of unity we can achieve that
$$ \zeta_i \in \QQ(e^{2\pi i/N}) \quad \text{where} \quad
N = \prod_{p \leq 2kg\,\, \text{prime}} p^{\nu_0(p)} \quad
\text{and} \quad \nu_0(p) = \lceil\log_p(1+\frac{g}{p-1}))\rceil.$$
\newline
In particular the $\zeta_i$ appearing in such a relation with 
the normalization $\zeta_1 =1$ belong to a finite set.
\end{Lemma} 
\par
{\bf Proof:} The following argument extends a theorem of Mann (\cite{Mn65})
from the case of rational coefficients to the case of coefficients
in a field of bounded degree over $\QQ$.
\newline
Suppose the irreducible relation has $\zeta_i \in \QQ(e^{2\pi i/N})$
and $N=p^\nu N'$ for some $\nu \geq \nu_0 +1$ and ${\rm gcd}(p,N')=1$.
Let $\zeta$ be a primitive $p^\nu$-th root of unity and let $\rho$ be a
primitive $N/p$-th root of unity. Resorting the relation according
to powers of $\zeta$ we obtain
\begin{equation}
 \sum_{j=0}^{p-1} b_j \zeta^j = 0 \quad \text{where} \quad
b_j = \sum_{i \in \Lambda_j} a_i \rho^{\alpha_i} 
\end{equation}
where $\Lambda_j=\{i: \exists \alpha_i \in \NN 
\,\text{such that}\, \rho^{\alpha_i} \zeta^j = \zeta_i \}$. 
The coefficients $b_j$ belong to $L=K(\rho)$. Since $K \subset
\QQ^{\rm ab}$ and $[K:\QQ] \leq g$ we know that 
$$K \subset \QQ(e^{2\pi i /p^{j_0(p)}}, p \,\; {\rm prime}).$$
Since cyclotomic fields for powers of different primes 
are linearly disjoint over $\QQ$ we deduce that 
$$ [L(\zeta):L] = [\QQ(\rho,\zeta), \QQ(\rho)] = p.$$
Hence $b_j=0$ for $j=0,\ldots,p-1$. Since the original 
relation was irreducible, this is only possible if
$\Lambda_{j_0} = \{1,\ldots,k\}$ for some $j_0$ (and the
other $\Lambda_j$ are empty). This means that
we can reduce $N$ by multiplying the original relation with
a suitable power of $\zeta$. 
\par
We have bounded the exponents that occur in the factorization
of $N$. It remains to
bound the size of primes dividing $N$. Suppose that
$p$ is prime and divides $N$ to the order
$\nu \leq \nu_0(p)$. As above let $\zeta$ be a primitive
$p^\nu$-th root of unity, but we let now $\rho$ be a primitive
$N/p^\nu$-th root of unity. Resorting the relation according
to powers of $\zeta$ we obtain
\begin{equation}
f(\zeta) := \sum_{j=0}^{p^\nu-1} b_j \zeta^j = 0 \quad \text{where} \quad
b_j = \sum_{i \in \Lambda_j} a_i \rho^{\alpha_i} 
\end{equation}
The coefficients $b_j$ of $f$ lie in $K(\rho)$. Since $\QQ(\rho) \cap
\QQ(\zeta) =\QQ$ the polynomial $f$ is a multiple of
the minimal polynomial $f_{\zeta/K}$ of $\zeta$ over $K$,
which has degree at least $\phi(p^\nu)-g$. Here $\phi$ denotes
Euler's $\phi$-function. On the other
hand by construction at most $k$ of the coefficients $b_j$ 
are non-zero. Hence there is somewhere a gap of size $p^\nu/k$
between non-zero $b_j$. Multiplying the relation by a suitable
power $\zeta$ we may suppose from the beginning that
$$\deg f \leq p^\nu(1-\frac{1}{k})$$
If $p^\nu/k -1  \geq g + p^{\nu-1}$ this leads to a contradiction
to the degree of $f_{\zeta/K}$. This condition is fulfilled
if the rough bound $p \leq 2kg$ is violated.
\hfill $\Box$
\par
{\bf Proof of Prop.\ \ref{Nbounded}:} Suppose the
finiteness statement was wrong. Then there exists a sequence
$(b^h_{i,n}, x_{i,n})$ for $n\in \NN$
satisfying $(Eq:e'')$ for all $e=1,\ldots,g-1$ and such that least
common multiple $N(n)$ of the orders of the $x_{i,n}$
is unbounded. We interpret the (solutions of the) equations
as relations between roots of unity
$$ \sum_{i=1}^g b^h_{i,n} x^e_{i,n} + \sum_{i=-g}^{-1}
b^h_{i,n} x^e_{i,n}= 0$$
with the convention that $b^h_{i,n} = b^h_{-i,n}$
and $x_{i,n}= x_{-i,n}^{-1}$.
\newline
For each $n$ and each $e$ we may write the relation in a (non-unique)
way as a sum of irreducible relations. The summands occuring
in such an irreducible relation form a partition of $I:= \{-g,\ldots,-1,
1,\ldots,g\}$. Since this set admits only finitely many partitions
we pass to a subsequence of $(b^h_{i,n}, x_{i,n})$
and suppose without loss of generality
that there are partitions $P_e$ consisting of subsets $P_{e,j}$
of $I$ such that for $e=1,\ldots,g-1$, for all $j$ and for
all $n \in \NN$ 
$$ \sum_{i \in P_{e,j}}  b^h_{i,n} x^e_{i,n}=0$$
is an irreducible relation. We apply Lemma \ref{MannThm} to these
relations and write 
\begin{equation} \label{useMann}
x_{i,n}^e = \zeta_{i,e,n} \sigma_{i,e,n} 
\end{equation}
with the following two properties: First, the
$\zeta_{i,e,n}$ are roots of unity of order bounded
by a function depending only on $g$ since the relations
are of length $\leq 2g$. Second if $i$ and $i'$ are both in $P_{e,j}$ then 
$\sigma_{i,j,n} = \sigma_{i',j,n}$. Passing to a subsequence
again we may suppose 
$$ \zeta_{i,e,n} = \zeta_{i,e} \quad \text{for all} \quad n \in \NN. $$
\newline
We want to limit the possible choices for the $\sigma_{i,e,n}$
to a finite set in order to obtain a contradiction.
From (\ref{useMann}) we deduce that
$$ \sigma_{i,e,n} = \sigma^e_{i,1,n} \frac{\zeta^e_{i,1}}{\zeta_{i,e}}.$$
This means that the $\sigma_{i,e,n}$ for different second arguments
are closely related. Since they coincide when the first argument
varies in a fixed partition set $P_{e,j}$ there is, roughly
speaking, a partition of $I$ that controls the $\sigma_{i,e,n}$
for all $e$ simultaneously. More precisely, consider the
following equivalence relation: 
$i\sim i'$ if there exists $(e,j)$ 
such that $P_{e,j} \supset \{i,i'\}$. Denote the corresponding
partition by $P_0 = \cup_j \,P_{0,j}$. Then 
$\sigma_{i,e,n}$ and $\sigma_{i',e,n}$ differ 
for $P_{0,j} \supset \{i,i'\}$ at worst 
by a product of $\zeta_{i,e}$ and a $(g-1)!$-th root of unity.
We may suppose that they actually coincide by recording
the discrepancy in modified $\zeta_{i,e}$. I.e.\ we write
\begin{equation} \label{usePart}
x_{i,n}^e = \widetilde{\zeta_{i,e}} \sigma^e_{i,n} 
\end{equation}
with the following two properties: First,
$\widetilde{\zeta_{i,e}}$ is root of unity of order bounded
by a function depending only on $g$ and second if 
$i$ and $i'$ are both in $P_{0,j}$ then 
$\sigma_{i,n} = \sigma_{i',n}$.
\newline
Suppose $P_{0,j}$ does not contain a pair $\{i,-i\}$. Then
the cardinality $k$ of $P_{0,j}$ is at most $g$. We fix $n \in \NN$.
The $b_{i,n}^h$ are solutions of the system of
equations
\begin{equation}
\sum_{i \in P_{o,j}} b_{i,n}^h x_{i,n}^{e} = 0, \quad e \in \{1,\ldots,k-1\}.
\end{equation}
Since the $b_{i,n}^h$ are real we may take complex conjugates to
see that they also solve the system of
equations for $ e \in \{-k+1,\ldots,-1\}$.
Since $x_{i,n} \neq x_{j,n}$ for $i \neq j$ the only solution
is $\sum_{i \in P_{o,j}} x_{i,n}^e = 0$, i.e.\ all $b_{i,n}^h$ for $i \in P_{0,j}$
are equal. This contradicts that the $b^h_i$ form a $\QQ$-basis
of $K$.
\par
On the other hand if $P_{0,j}$ contains $\{i_0, -i_0\}$ we deduce from
$$ \zeta_{-i_0,1} \sigma_{i_0,n} =
\zeta_{-i_0,1} \sigma_{-i_0,n} = x_{-i_0,n} = x_{i_0,n}^{-1} = 
(\zeta_{i_0,1}  \sigma_{i_0,n})^{-1} $$
that $\sigma_{i_0,n}$ runs for $n \in \NN$
through a finite set. By construction the same holds for
$\sigma_{i,n}$ for all $i \in P_{0,j}$.
Applying this to all partition sets of $P_0$ the orders of the
$\sigma_{i,n}$ and hence the orders of the $x_{i,n}$ are bounded.
This contradicts the choice of the sequence $x_{i,n}$.
\hfill $\Box$
\par


\par
Martin M{\"o}ller: Universit{\"a}t Essen, FB 6 (Mathematik) \newline 
45117 Essen, Germany \newline
e-mail: martin.moeller@uni-essen.de \newline

\end{document}